\documentclass{article}

\usepackage{amsmath,amssymb}

\begin{document}

\title{Rates of convergence for nearest neighbor estimators with the smoother regression function}

\author{Takanori Ayano\footnote{t-ayano@cr.math.sci.osaka-u.ac.jp}\\Department of Mathematics, Osaka University}

\date{}

\maketitle

\begin{abstract}
Let $(X,Y)$ be a $\mathbb{R}^d\times\mathbb{R}$-valued random vector. 
In regression analysis one wants to estimate the regression function $m(x):={\bf E}(Y|X=x)$ from a data.    
In this paper we consider the rate of convergence for the $k$ nearest neighbor estimator 
in case that $X$ is uniformly distributed on $[0,1]^d$, $\mbox{{\bf Var}}(Y|X=x)$ is bounded, and $m$ is $(p,C)$-smooth. 
It is an open problem whether the optimal rate can be achieved by some $k$ nearest neighbor estimator in case of $1<p\le1.5$. 
We solve the problem affirmatively. This is the main result of this paper. 
Throughout this paper, we assume that the data is independent and identically distributed and as an error criterion we use the expected $L_2$ error.
\end{abstract}

\vspace{1ex}


\section{Introduction}
Let $(X,Y)$ be a $\mathbb{R}^d\times\mathbb{R}$-valued random vector. 
In regression analysis, one wants to predict the value of $Y$ after having observed the value of $X$, i.e.
to find a measurable function $f$ such that the mean squared error ${\bf{E}}_{XY}\left(f(X)-Y\right)^2$ is minimized, 
where ${\bf E}_{XY}$ denotes the expectation with respect to $(X,Y)$.
Let $m(x):={\bf E}\{Y|X=x\}$ (regression function), which is the conditional expectation of $Y$ given $X=x$.  
Then $m(x)$ is the solution of the minimization problem. 
In fact, one can check for any measurable function $f$, 
\[{\bf E}_{XY}\left(f(X)-Y\right)^2={\bf E}_{XY}\left(m(X)-Y\right)^2+{\bf E}_X\left(f(X)-m(X)\right)^2.\]
In statistics, only the data is available, (the distribution of $(X,Y)$ and $m$ are not available), and one needs to estimate the function $m$ from the data $\{(X_i,Y_i)\}_{i=1}^n$, which 
are independently distributed according to the distribution of $(X,Y)$. 
We wish to construct an estimator $m_n$ of $m$ such that the expected $L_2$ error $R(m_n):={\bf E}_{X^nY^n}{\bf E}_X\left(m_n(X)-m(X)\right)^2$
 is as small as possible, where ${\bf E}_{X^nY^n}$ denotes the expectation with respect to the data.
In order to analyze the performance of estimators theoretically, it is very important to evaluate how fast the error $R(m_n)$ converges to zero, when the data size $n$ tends to infinity. 
In this paper we consider $k$-NN (nearest neighbor) estimators and the rate of convergence in case that $m$ is $(p,C)$-smooth (cf.Gy\"orfi et al., 2002, p.37). 

The $k$-NN estimator is defined as follows.
Given $x\in\mathbb{R}^d$, we rearrange the data $(X_1,Y_1),\dots,(X_n,Y_n)$ in the ascending order of the values of $\|X_i-x\|$. 
As a tie-breaking rule, if $\|X_i-x\|=\|X_j-x\|$ and $i<j$, we declare that $X_i$ is ``closer" to $x$ than $X_j$. 
We write the rearrange sequence by $\left(X_{1,x},Y_{1,x}\right),\dots,\left(X_{n,x},Y_{n,x}\right)$. 
Notice that $\{(X_{i,x},Y_{i,x})\}_{i=1}^n$ is expressed by $\{(X_{\pi(i)},Y_{\pi(i)})\}_{i=1}^n$ using a permutation $\pi:\{1,\dots,n\}\to\{1,\dots,n\}$ depending on $x\in\mathbb{R}^d$. 
Then for $1\le k\le n$, the $k$-NN estimator $m_n$ is defined by 
\[m_n(x)=\frac{1}{k}\sum_{i=1}^{k}Y_{i,x}\]
For the details about $k$-NN estimators, for example, see Chapter 6 in Gy\"orfi et al. (2002). 

Let $p,C>0$, and express $p$ by $p=q+r,\;q\in{\mathbb{Z}_{\ge 0}},\;0<r\le 1$. We say that a function $m:\mathbb{R}^d\to\mathbb{R}$ is $(p,C)$-smooth
if for all $q_1,\dots,q_d\in\mathbb{Z}_{\ge 0}$ with $q=q_1+\cdots+q_d$, the partial derivatives
$\frac{\partial^qm}{\partial x_1^{q_1}\cdots\partial x_d^{q_d}}$
exist and for all $x,z\in\mathbb{R}^d$ the following is satisfied. 
\[\left|\frac{\partial^qm}{\partial x_1^{q_1}\cdots\partial x_d^{q_d}}(x)-\frac{\partial^qm}{\partial x_1^{q_1}\cdots\partial x_d^{q_d}}(z)\right|\le C\|x-z\|^r\]

\vspace{1ex}

For $p,C,\sigma>0$, let $\mathcal{D}(p,C,\sigma)$ be the class of distributions of $(X,Y)$ such that:  

\vspace{1ex}

(I) $X$ is uniformly distributed on $[0,1]^d$; 

\vspace{1ex}

(II) $\mbox{{\bf Var}}(Y|X=x)\le\sigma^2$; 

\vspace{1ex}

(III) $m$ is ($p,C$)-smooth, 

\vspace{1ex}

where $\mbox{{\bf Var}}(Y|X=x)$ denotes the variance of $Y$ given $X=x$. 

\vspace{1ex}

The lower bound for the class $\mathcal{D}(p,C,\sigma)$ is known (cf.Gy\"orfi et al., 2002, p.38):
\begin{equation}
\liminf_{n\to\infty}\inf_{m_n}\sup_{(X,Y)\in\mathcal{D}(p,C,\sigma)}n^{2p/(2p+d)}R(m_n)\ge\mbox{const.}>0\label{stone}
\end{equation}
where $\inf_{m_n}$ denotes the infimum over all the estimators.  

For $0<p\le1$, the rate $n^{-2p/(2p+d)}$ is achieved by the $k$-NN estimator (cf.Gy\"orfi et al., 2002, pp.93,99):  
\[\sup_{(X,Y)\in\mathcal{D}(p,C,\sigma)}R(m_n)\le \mbox{const.}\;n^{-2p/(2p+d)}\]
For $p>1.5$, it is shown that the rate $n^{-2p/(2p+d)}$ is unachievable by any $k$-NN estimator  
and it is presented as a conjecture that even for $1<p\le1.5$, the rate $n^{-2p/(2p+d)}$ will be achieved by some $k$-NN estimator (cf.Gy\"orfi et al., 2002, p.96). 
In this paper, we show that the conjecture is right (Theorem). 
Regression analysis is used in many fields for example economics, medicine, pattern recognition  etc. (cf.Gy\"orfi et al., 2002, pp.4-9).  
Nearest neighbor estimators are very important in regression analysis.  
We have shown the performance of the nearest neighbor estimator theoretically.

Throughout this paper we will use the following notations : $\mathbb{R},\mathbb{R}_{>0},\mathbb{Z}_{\ge0},\mathbb{N}$ are the sets of reals, positive reals, nonnegative integers and positive intgers. 
For a measurable set $D\subset\mathbb{R}^d$, $\mbox{vol}(D)$ denotes the Lebesgue measure of $D$. 
For $x\in\mathbb{R}^d$, $\|x\|$ denotes the Euclidean norm of $x$. For $u,v\in\mathbb{R}^d$, we define $H(u,v):=\{w\in\mathbb{R}^d\;|\;\|w-u\|\le\|v-u\|\}$ and $G(u,v):=H(u,v)\cap[0,1]^d$.  
For $a>0$, $\lfloor a\rfloor$ denotes $b$ such that $b\le a<b+1$.

\section{Related Work}
In this section, we overview  the related work about consistency and the rate of convergence. 
For consistency, it was shown in Stone (1977) that
the $k$-NN estimators are universally consistent. Since
then it was shown that many estimators share this property (cf.Devroye et al.,1994, Greblicki et al.,1984, Gy\"orfi and Walk, 1997, Kohler, 1999, Kohler and Krzy\.zak, 2001, Kohler, 2002, Lugosi and Zeger, 1995, Nobel, 1996, Walk, 2002, Walk, 2005, Walk, 2008). 
For the rate of convergence, we know several results as follow: 
\begin{itemize}
\item Stone (1982) proved the lower bound (1);
\item for the distributions
satisfying (II)(III) with $0<p\leq 1$ and the partitioning, kernel, and $k$-NN estimator,
 the rate $n^{-2p/(2p+d)}$ is achievable if $X$ is bounded (for the $k$-NN estimator, the condition $d > 2p$ is required as well) (cf.Gy\"orfi, 1981, Gy\"orfi et al., 2002, Kulkarni and Posner, 1995, Spiegelman and Sacks, 1980); 
\item Kohler et al.(2006, 2009) proved the same statement later without assuming that $X$ should be bounded;
\item  for the partitioning estimators and 
the class ${\mathcal D}(p,C,\sigma)$ with $p > 1$,  
the rate $n^{-2p/(2p+d)}$ is unachievable (cf.Gy\"orfi et al., 2002);
\item for the kernel estimators, the rate $n^{-2p/(2p+d)}$ is achievable for ${\mathcal D}(p,C,\sigma)$ with $0<p\leq 1.5$ and
 is unachievable for that with $p>1.5$ (cf.Gy\"orfi et al., 2002);
\end{itemize}

If we summarize the above results in Table 1, only the following problem remains:
Does the $k$-NN estimator achieve the rate $n^{-2p/(2p+d)}$ under (II)(III) even for $1<p\leq 1.5$ ?
The problem is still hard, but we solve the statement affirmatively under (I)(II)(III).

\vspace{1ex}

Table 1 : the achievability of $n^{-2p/(2p+d)}$ for the estimators and ${\mathcal D}(p,C,\sigma)$
\begin{center}
\begin{tabular}{lcc}\hline
{} & achievable & unachievable \\ \hline
partitioning & $0<p\le 1$  & $p>1$  \\ \hline
kernel & $0<p\le1.5$ & $p>1.5$ \\ \hline
$k$-NN & $0<p\le1$  & $p>1.5$ \\ \hline
\end{tabular}
\end{center}

\section{Main Result}
For $p,C,\sigma>0$, let $\mathcal{D}(p,C,\sigma)$ be the class of distributions of $(X,Y)$ such that:  

\vspace{1ex}

(I) $X$ is uniformly distributed on $[0,1]^d$; 

\vspace{1ex}

(II) $\mbox{{\bf Var}}(Y|X=x)\le\sigma^2$; 

\vspace{1ex}

(III) $m$ is ($p,C$)-smooth, 

\vspace{1ex}

where $\mbox{{\bf Var}}(Y|X=x)$ denotes the variance of $Y$ given $X=x$.

\vspace{1ex}

Then we get the following theorem: 

\vspace{1ex}

{\bf{Theorem}}

Let $1<p\le 1.5$ and let $m_n$ be the $k$-NN estimator with $k=\lfloor n^{2p/(2p+d)}\rfloor$. 
Then there exists $C_1>0$ (which does not depend on $n$) such that 
\[\sup_{(X,Y)\in\mathcal{D}(p,C,\sigma)}{\bf E}_{X^nY^n}{\bf E}_X\left(m_n(X)-m(X)\right)^2\le C_1n^{-2p/(2p+d)}.\]

\section{Proof of Theorem}
Suppose we are given $X=x,X_1=x_1,\dots,X_n=x_n$. 
We take the expectation with respect to $Y_1,\dots,Y_n$. Then the following bias-variance decomposition is well-known (cf.Gy\"orfi et al., 2002, p.94): 
\[{\bf E}_{Y^n}\left(m_n(x)-m(x)\right)^2={\bf E}_{Y^n}\left(\frac{1}{k}\sum_{i=1}^k\left(Y_{i,x}-m(x)\right)\right)^2\]
\[={\bf E}_{Y^n}\left(\frac{1}{k}\sum_{i=1}^k\left(Y_{i,x}-m(x_{i,x})\right)\right)^2+\left\{\frac{1}{k}\sum_{i=1}^k\left(m(x_{i,x})-m(x)\right)\right\}^2\] 
\begin{equation}
\le \frac{\sigma^2}{k}+\frac{1}{k^2}\left\{\sum_{i=1}^k\left(m(x_{i,x})-m(x)\right)\right\}^2. \;\;\;\;(\because \mbox{(II)})\hspace{14ex}\label{eq:fixxX}
\end{equation}

\vspace{1ex}

We evaluate the second term of (\ref{eq:fixxX}).  
Let $x_{i,x}=(x_{i,x}^{(1)},\dots,x_{i,x}^{(d)})$ and $x=(x^{(1)},\dots,x^{(d)})$. Let $m_s$ be the partial derivative of $m$ with respect to the $s$-th component. 
Then by the mean-value theorem, there exists $u_i\in\mathbb{R}^d$ such that $\|u_i-x\|\le \|x_{i,x}-x\|$ and 
\[\left\{\sum_{i=1}^k(m(x_{i,x})-m(x))\right\}^2=\left\{\sum_{i=1}^k\sum_{s=1}^dm_s(u_i)(x_{i,x}^{(s)}-x^{(s)})\right\}^2,\]
(the idea using the mean-value theorem is due to Gy\"orfi et al., 2002, p.84) by Cauchy-Schwarz's inequality
\[\le2\left\{\sum_{i=1}^k\sum_{s=1}^d(m_s(u_i)\!-\!m_s(x))(x_{i,x}^{(s)}\!-\!x^{(s)})\right\}^2\!+\!2\left\{\sum_{s=1}^dm_s(x)\sum_{i=1}^k(x_{i,x}^{(s)}\!-\!x^{(s)})\right\}^2\] 
\[\le2kd\sum_{i=1}^k\sum_{s=1}^d(m_s(u_i)\!-\!m_s(x))^2(x_{i,x}^{(s)}\!-\!x^{(s)})^2\!+\!2d\sum_{s=1}^dm_s(x)^2\left\{\sum_{i=1}^k(x_{i,x}^{(s)}\!-\!x^{(s)})\right\}^2,\]
let $L>0$ such that $\max_{1\le s\le d,x\in [0,1]^d}|m_s(x)|\le L$, because $m$ is $(p,C)$-smooth and $\|u_i-x\|\le \|x_{i,x}-x\|$, 
\[\le2kdC^2\sum_{i=1}^k\sum_{s=1}^d\|x_{i,x}-x\|^{2p-2}(x_{i,x}^{(s)}-x^{(s)})^2+2dL^2\sum_{s=1}^d\left\{\sum_{i=1}^k(x_{i,x}^{(s)}-x^{(s)})\right\}^2\]
\[=2kdC^2\sum_{i=1}^k\|x_{i,x}-x\|^{2p}+2dL^2\sum_{i=1}^k\|x_{i,x}-x\|^2\hspace{26ex}\]
\[+2dL^2\sum_{s=1}^d\sum_{1\le i\neq j\le k}(x_{i,x}^{(s)}-x^{(s)})(x_{j,x}^{(s)}-x^{(s)})\hspace{24ex}\]
We regard $x,x_1,\dots,x_n$ as the random variables $X,X_1,\dots,X_n$ and take the expectation with respect to $X,X_1,\dots,X_n$. 
\[{\bf E}_X{\bf E}_{X^nY^n}\left(m_n(X)-m(X)\right)^2\]
\begin{equation}
\le \frac{\sigma^2}{k}+\frac{2dC^2}{k}{\bf E}_X{\bf E}_{X^n}\sum_{i=1}^k\|X_{i,X}-X\|^{2p}+\frac{2dL^2}{k^2}{\bf E}_X{\bf E}_{X^n}\sum_{i=1}^k\|X_{i,X}-X\|^2\label{eq:nofix1}
\end{equation}
\begin{equation}
+\frac{2dL^2}{k^2}{\bf E}_X{\bf E}_{X^n}\sum_{s=1}^d\sum_{1\le i\neq j\le k}\left(X_{i,X}^{(s)}-X^{(s)}\right)\left(X_{j,X}^{(s)}-X^{(s)}\right)\label{eq:nofix2}
\end{equation}

\vspace{1ex}

In order to evaluate the second and third terms in (\ref{eq:nofix1}), the following proposition is available. 

\vspace{1ex}

{\bf{Proposition}} (Gy\"orfi et al., 2002, pp.95,99)

For any $\gamma>0$, there exists $c_1>0$ (depending on $\gamma$ and $d$) such that,  
\[\frac{1}{k}{\bf E}_X{\bf E}_{X^n}\sum_{i=1}^k\|X_{i,X}-X\|^{2\gamma}\le c_1\left(\frac{k}{n}\right)^{2\gamma/d}.\]

The proposition is proved originally for $\gamma=1$ in Gy\"orfi et al., 2002, but we have extended it to the general $\gamma>0$. 
We proceed to evaluate (\ref{eq:nofix2}).   

Let $D=\{(x,x_1,\dots,x_n)\;|\;\|x_i-x\|<\|x_{k+1}-x\|,i=1,\dots, k,\|x_j-x\|>\|x_{k+1}-x\|,j=k+2,\dots, n\}$. 

\vspace{1ex}

{\bf{Claim 1}}
\[{\bf E}_X{\bf E}_{X^n}\sum_{1\le i\neq j\le k}\left(X_{i,X}^{(s)}-X^{(s)}\right)\left(X_{j,X}^{(s)}-X^{(s)}\right)\]
\[=\frac{n\cdots(n-k)}{k!}\int_D\sum_{1\le i\neq j\le k}\left(x_i^{(s)}-x^{(s)}\right)\left(x_j^{(s)}-x^{(s)}\right)dx_1\cdots dx_ndx\]
(See Appendix for proof) 

\vspace{2ex}

From Claim 1, 

since $x_1,\dots,x_k\in G(x,x_{k+1})$ and $x_{k+2},\dots,x_n\in [0,1]^d\backslash G(x,x_{k+1})$ on $D$, 
\[T:=\frac{1}{k^2}{\bf E}_X{\bf E}_{X^n}\sum_{s=1}^d\sum_{1\le i\neq j\le k}\left(X_{i,X}^{(s)}-X^{(s)}\right)\left(X_{j,X}^{(s)}-X^{(s)}\right)\]
\[=\frac{n\cdots(n-k)}{k^2k!}\sum_{s=1}^d\sum_{1\le i\neq j\le k}\int_D\left(x_i^{(s)}-x^{(s)}\right)\left(x_j^{(s)}-x^{(s)}\right)dx_1\cdots dx_ndx,\]
\[=\frac{n\cdots(n-k)}{k^2k!}\sum_{s=1}^d\sum_{1\le i\neq j\le k}\int_{[0,1]^d}dx\int_{[0,1]^d}dx_{k+1}\int_{G(x,x_{k+1})}(x_i^{(s)}-x^{(s)})dx_i\]
\[\int_{G(x,x_{k+1})}(x_j^{(s)}-x^{(s)})dx_j\cdot\mbox{vol}[G(x,x_{k+1})]^{k-2}(1-\mbox{vol}[G(x,x_{k+1})])^{n-k-1}.\]

\vspace{1ex}

Let $W:=\{(x,x_{k+1})\;|\;G(x,x_{k+1})\neq H(x,x_{k+1})\}$. Since for $(x,x_{k+1})\notin W$, 
$\int_{G(x,x_{k+1})}(x_i^{(s)}-x^{(s)})dx_i=\int_{H(x,x_{k+1})}(x_i^{(s)}-x^{(s)})dx_i=0$, we obtain 
\[T=\sum_{s=1}^d\frac{n\cdots(n-k)}{k^2(k-2)!}\int_Wdx_{k+1}dx\left(\int_{G(x,x_{k+1})}(x_1^{(s)}-x^{(s)})dx_1\right)^2\]
\[\mbox{vol}[G(x,x_{k+1})]^{k-2}(1-\mbox{vol}[G(x,x_{k+1})])^{n-k-1}.\]

\vspace{1ex}

{\bf{Claim 2}}\hspace{8ex}There exists $c_2>0$ (depending only on $d$) such that 
\[\left|\int_{G(x,x_{k+1})}(x_1^{(s)}-x^{(s)})\;dx_1\right|\le c_2\cdot\mbox{vol}[G(x,x_{k+1})]^{(d+1)/d}\]
(See Appendix for proof)

\vspace{2ex}
 
From Claim 2, 
\[T\!\le\! c_2^2d\frac{n\cdots(n\!-\!k)}{k^2(k\!-\!2)!}\!\int_W\mbox{vol}[G(x,x_{k+1})]^{k+\frac{2}{d}}\left\{1\!-\!\mbox{vol}[G(x,x_{k+1})]\right\}^{n-k-1}dx_{k+1}dx\]
\
\begin{equation}
=c_2^2d\frac{n\cdots(n-k)}{k^2(k-2)!}\int_0^1u^{k+\frac{2}{d}}(1-u)^{n-k-1}dF(u)\hspace{24ex}\label{density}
\end{equation}
where $F(u)$ is the Lebesgue measure of $S(u):=\{(x,x_{k+1})\in W\;|\;0\le\mbox{vol}[G(x,x_{k+1})]\le u\}$ for $0\le u\le 1$. 

\vspace{1ex}

{\bf{Claim 3}}
\[F(u)\!=\!\left\{\begin{array}{ll}\displaystyle{\!2d\!\sum_{i=0}^{d-1}\!\left\{\frac{{}_{d-1}C_i(-2)^{d-1-i}}{(d\!-\!i)e_2^{(d-i)/d}}\!-\!\frac{{}_{d-1}C_i(-2)^{d-1-i}}{(2d\!-\!i)e_2^{(d-i)/d}}\right\}\!u^{(2d-i)/d}} & \!\left(0\le u\le \displaystyle{\frac{e_2}{2^d}}\right) \\ \displaystyle{u-2de_2\int_0^{1/2}(x^{(1)})^d(1-2x^{(1)})^{d-1}dx^{(1)}} & \left(\displaystyle{\frac{e_2}{2^d}}\le u\le 1\right)\end{array}\right.\]
where $e_2=\mbox{vol}[\{y\in\mathbb{R}^d\;|\;\|y\|\le1\}]$, and $e_2\le2^d$. 

\vspace{1ex}

(See Appendix for proof)

\vspace{2ex}

For $0<u<e_2/2^d$ and $e_2/2^d<u<1$, let $f(u):=F'(u)$. ($f(u)\ge0$)   
\[f(u)\!=\!\left\{\!\begin{array}{ll}\displaystyle{\sum_{i=0}^{d\!-\!1}(4d\!-\!2i)\!\left\{\frac{{}_{d-1}C_i(\!-\!2)^{d\!-\!1\!-\!i}}{(d\!-\!i)e_2^{(d\!-\!i)/d}}\!-\!\frac{{}_{d\!-\!1}C_i\;(\!-\!2)^{d\!-\!1\!-\!i}}{(2d\!-\!i)e_2^{(d\!-\!i)/d}}\right\}\!u^{(d\!-\!i)/d}} & \!\left(0\!<\!u\!<\!\displaystyle{\frac{e_2}{2^d}}\right) \\ 1 & \!\left(\displaystyle{\frac{e_2}{2^d}}\!<\!u\!<\!1\right)\end{array}\right.\]
Let $f(0)=f(e_2/2^d)=f(1)=0$.  
There exists $c_3>0$ (depending only on $d$) such that $f(u)\le c_3u^{1/d}$, because, 
for $e_2/2^d<u<1$, $f(u)=1\le (2/e_2^{1/d})u^{1/d}$ and for the other $u$ it is trivial. 

For $\alpha\in\mathbb{R}_{>0}$ and $\beta\in\mathbb{N}$, let $B(\alpha,\beta):=\displaystyle{\int_0^1u^{\alpha-1}(1-u)^{\beta-1}du}$ (Beta function).
Then the following formula is well-known: 
\[B(\alpha,\beta)=\frac{\Gamma(\alpha)\Gamma(\beta)}{\Gamma(\alpha+\beta)}=\frac{\Gamma(\alpha)\;(\beta-1)!}{(\alpha+\beta-1)\cdots\alpha\Gamma(\alpha)}=\frac{(\beta-1)!}{\alpha\cdots(\alpha+\beta-1)},\]
where $\Gamma$ is Gamma function. 

\vspace{1ex}

On the other hand, 
\[\lim_{n\to\infty}\frac{n!}{(1+\frac{3}{d})\cdots(n+\frac{3}{d})}\cdot n^{\frac{3}{d}}=\lim_{n\to\infty}\frac{\Gamma(n+1)\cdot\Gamma(1+\frac{3}{d})}{\Gamma(n+1+\frac{3}{d})}\cdot n^{\frac{3}{d}}\] 
By Stirling's formula,
\[=\lim_{n\to\infty}\Gamma(1+\frac{3}{d})\frac{\sqrt{2\pi n}\left(\frac{n}{e}\right)^n}{\sqrt{2\pi(n+\frac{3}{d})}\left(\frac{n+\frac{3}{d}}{e}\right)^{n+\frac{3}{d}}}\cdot n^{\frac{3}{d}}=\Gamma(1+\frac{3}{d})\]
Therefore, there exist $c_4,c_5>0$ (depending only on $d$) such that 
\[c_4n^{-\frac{3}{d}}\le\frac{n!}{(1+\frac{3}{d})\cdots(n+\frac{3}{d})}\le c_5n^{-\frac{3}{d}}\]
From (\ref{density}), 
\[T\le c_2^2d\frac{n\cdots(n-k)}{k^2(k-2)!}\int_0^1u^{k+\frac{2}{d}}(1-u)^{n-k-1}f(u)du\hspace{26ex}\]
\[\le c_2^2c_3d\frac{n\cdots(n-k)}{k^2(k-2)!}\int_0^1u^{k+\frac{3}{d}}(1-u)^{n-k-1}du\hspace{26ex}\]
\[=c_2^2c_3d\frac{n\cdots(n-k)}{k^2(k-2)!}B\left(k+1+\frac{3}{d},\;n-k\right)\hspace{29ex}\]
\[=c_2^2c_3d\frac{n\cdots(n-k)}{k^2(k-2)!}\frac{(n-k-1)!}{(k+1+\frac{3}{d})\cdots(n+\frac{3}{d})}\le c_2^2c_3d\frac{(k+1)\cdots n}{(k+1+\frac{3}{d})\cdots(n+\frac{3}{d})}\]
\begin{equation}
=c_2^2c_3d\frac{n!}{(1+\frac{3}{d})\cdots(n+\frac{3}{d})}/\frac{k!}{(1+\frac{3}{d})\cdots(k+\frac{3}{d})}\le\frac{c_2^2c_3c_5d}{c_4}\left(\frac{k}{n}\right)^{3/d}\hspace{6ex}\label{ay}
\end{equation}

\vspace{1ex}

Therefore, from (\ref{eq:nofix1}), (\ref{eq:nofix2}), (\ref{ay}), and Proposition, there exist $C_2,C_3,C_4>0$ (which do not depend on $n$) such that  
\[{\bf E}_X{\bf E}_{X^nY^n}\left(m_n(X)-m(X)\right)^2\le \frac{\sigma^2}{k}+C_2\left(\frac{k}{n}\right)^{2p/d}+\frac{C_3}{k}\left(\frac{k}{n}\right)^{2/d}+C_4\left(\frac{k}{n}\right)^{3/d}\]
Assuming $p\le 1.5$, if we set $k=\lfloor n^{2p/(2p+d)}\rfloor$, there exists $C_1>0$ (which does not depend on $n$) such that  
\[{\bf E}_X{\bf E}_{X^nY^n}\left(m_n(X)-m(X)\right)^2\le C_1n^{-2p/(2p+d)}\]
We have got Theorem. 

\hspace{74ex}$\square$



\vspace{4ex}

{\bf \Large Appendix}

\appendix

\section{Proof of Claim 1}
Let $h\in N:=\{1,\dots,n\}$ and $I,J\subset N\backslash\{h\}$ such that $\sharp I=k, I\cap J=\{\}\mbox{(empty)}, I\cup J=N\backslash\{h\}$, 
where, $\sharp$ denotes the number of the elements.  
Let $D(I,J,h):=\{(x,x_1,\dots,x_n)\;|\;\|x_i-x\|<\|x_h-x\|,i\in I,\|x_j-x\|>\|x_h-x\|,j\in J\}$.
Since $\mbox{vol}\left\{[0,1]^{d(n+1)}\backslash\cup_{I,J,h}D(I,J,h)\right\}=0$ and for $(I,J,h)\neq(I',J',h'), D(I,J,h)\cap D(I',J',h')=\{\}$, we have 
\[{\bf E}_X{\bf E}_{X^n}\sum_{1\le i\neq j\le k}\left(X_{i,X}^{(s)}-X^{(s)}\right)\left(X_{j,X}^{(s)}-X^{(s)}\right)\]
\[=\sum_{I,J,h}\int_{D(I,J,h)}\sum_{i,j\in I,i\neq j}\left(x_i^{(s)}-x^{(s)}\right)\left(x_{j}^{(s)}-x^{(s)}\right)dx_1\cdots dx_ndx\]
Since for each $(I,J,h)$ the above integral has the same value and the number of $(I,J,h)$ is ${}_nC_k\cdot (n-k)$, we get Claim 1. 

\hspace{74ex}$\square$

\section{Proof of Claim 2}
Let $e_1:=\displaystyle{\int_{\|y\|\le1,\;y^{(s)}\ge0}y^{(s)}dy}$ and $e_2:=\displaystyle{\int_{\|y\|\le 1}dy}$, then for any $R\ge0$,
\[\int_{\|y\|\le R,\;y^{(s)}\ge0}y^{(s)}dy=e_1R^{d+1},\]
and 
\begin{equation}
\int_{\|y\|\le R}dy=e_2R^d\label{e_2}, 
\end{equation}
thus we have 
\[\left|\int_{G(x,x_{k\!+\!1})}\!(x_1^{(s)}\!-\!x^{(s)})dx_1\right|\!\le\!\int_{H(x,x_{k\!+\!1})}\!|x_1^{(s)}\!-\!x^{(s)}|dx_1\!=\!2\!\int_{\|y\|\le \|x\!-\!x_{k+1}\|,y^{(s)}\!\ge \!0}\!y^{(s)}dy\]
\[=2e_1\|x-x_{k+1}\|^{d+1}=\left\{\frac{(2e_1)^{d/(d+1)}}{e_2}\mbox{vol}[H(x,x_{k+1})]\right\}^{(d+1)/d}\hspace{10ex}\]

\vspace{1ex}

If we prove the following lemma, the proof of Claim 2 is complete: 

\vspace{1ex}

{\bf Lemma}

There exists $e_3$ (depending only on $d$) such that for any $u,v\in[0,1]^d$, 
\[\mbox{vol}[G(u,v)]\ge e_3\mbox{vol}[H(u,v)]\]
(Proof of Lemma)

Suppose $\|u-v\|\le 1/2$. 
Let $I:=\{i\;|\;0\le u^{(i)}\le 1/2\}$ and $M:=\{w\;|\;\|w-u\|\le\|u-v\|,\;w^{(i)}\ge u^{(i)},\;i\in I,\;w^{(j)}\le u^{(j)},\;j\notin I\}$. 
Then, $M\subset G(u,v)$ and $\mbox{vol}[M]=2^{-d}\mbox{vol}[H(u,v)]$, thus, 
\begin{equation}
\mbox{vol}[G(u,v)]\ge 2^{-d}\;\mbox{vol}[H(u,v)]\label{claim2.1} 
\end{equation} 
Suppose $\|u-v\|>1/2$. Since $\|u-v\|\le \sqrt{d}$, we have $\mbox{vol}[H(u,v)]\le e_2d^{d/2}$. 
From (\ref{claim2.1}), for $z\in\mathbb{R}^d$ such that $\|u-z\|=1/2$, 
\[\mbox{vol}[G(u,v)]\ge \mbox{vol}[G(u,z)]\ge2^{-d}\;\mbox{vol}[H(u,z)]=2^{-d}\left(e_22^{-d}\right).\] 
\[\therefore\;\;\frac{\mbox{vol}[G(u,v)]}{\mbox{vol}[H(u,v)]}\ge 2^{-2d}d^{-d/2}\]
Let $e_3:=\min\{2^{-d},2^{-2d}d^{-d/2}\}$, then we get Lemma. 

\hspace{74ex}$\square$ 

\section{Proof of Claim 3}
Let
\[V_{i}:=\left\{x\in[0,1]^d\;|\;x^{(i)}\le \min\{x^{(j)},1-x^{(j)}\},j=1,\dots d\right\}\]
\[V_{i+d}:=\left\{x\in[0,1]^d\;|\;1-x^{(i)}\le \min\{x^{(j)},1-x^{(j)}\},j=1,\dots d\right\}.\] 
Since $\cup_{i=1}^{2d}V_i=[0,1]^d$ and for $i\neq j, \mbox{vol}[V_i\cap V_j]=0$, by Fubini's theorem, 
\[F(u)\!=\!\int_{S(u)}dx_{k\!+\!1}\;dx=\int_{[0,1]^d}\left\{\int_{S(u)}dx_{k\!+\!1}\right\}dx\!=\!\sum_{i=1}^{2d}\int_{V_i}\left\{\int_{S(u)}dx_{k\!+\!1}\right\}dx\]

\vspace{1ex}

Without loss of generality, we assume $x\in V_1$.  

Let $y:=(0,x^{(2)},\dots,x^{(d)})$.
Since $H(x,y)\subset[0,1]^d$, for $u<\mbox{vol}[H(x,y)]$, 
\[0\le \mbox{vol}[G(x,x_{k+1})]\le u\Longrightarrow G(x,x_{k+1})=H(x,x_{k+1})\Longrightarrow (x,x_{k+1})\notin W\]
\[\Longrightarrow\{x_{k+1}\;|\;(x,x_{k+1})\in S(u)\}=\{\}\hspace{32ex}\]
For $\mbox{vol}[H(x,y)]\le u$ and $z\in\mathbb{R}^d$ such that $\mbox{vol}[G(x,z)]=u$, we have 
\[\{x_{k+1}\;|\;(x,x_{k+1})\in S(u)\}=G(x,z)\backslash H(x,y)\]
and from (\ref{e_2}), $\mbox{vol}[G(x,z)\backslash H(x,y)]=u-e_2(x^{(1)})^d.$
\[U:=\int_{V_1}\left\{\int_{S(u)}dx_{k+1}\right\}dx=\int_{V_1}\max\{u-e_2(x^{(1)})^d,0\}\;dx\]
\[=\int_0^{\min\left\{(u/e_2)^{1/d},1/2\right\}}\left\{\left(u-e_2(x^{(1)})^d\right)\int_{x^{(1)}}^{1-x^{(1)}}dx^{(2)}\cdots\int_{x^{(1)}}^{1-x^{(1)}}dx^{(d)}\right\}dx^{(1)}\]
\[=\int_0^{\min\left\{(u/e_2)^{1/d},1/2\right\}}\left(u-e_2(x^{(1)})^d\right)(1-2x^{(1)})^{d-1}dx^{(1)}\hspace{20ex}\]

\vspace{1ex}

For $(u/e_2)^{1/d}\le1/2$, i.e. $0\le u\le e_2/2^d$, @
\[U=\int_0^{(u/e_2)^{1/d}}\left(u-e_2(x^{(1)})^d\right)(1-2x^{(1)})^{d-1}dx^{(1)}\]
\[=\int_0^{(u/e_2)^{1/d}}\left(u-e_2(x^{(1)})^d\right)\left\{\sum_{i=0}^{d-1}\;{}_{d-1}C_i\;1^i\;(-2x^{(1)})^{d-1-i}\right\}dx^{(1)}\]
\[=\sum_{i=0}^{d-1}\left\{\frac{{}_{d-1}C_i(-2)^{d-1-i}}{(d-i)e_2^{(d-i)/d}}-\frac{{}_{d-1}C_i\;(-2)^{d-1-i}}{(2d-i)\;e_2^{(d-i)/d}}\right\}u^{(2d-i)/d}\hspace{10ex}\]

\vspace{1ex}

For $1/2\le(u/e_2)^{1/d}$, i.e. $e_2/2^d\le u\le 1$, 
\[U=\int_0^{1/2}\left(u-e_2(x^{(1)})^d\right)\;(1-2x^{(1)})^{d-1}dx^{(1)}\]
\[=\frac{u}{2d}-e_2\int_0^{1/2}(x^{(1)})^d\;(1-2x^{(1)})^{d-1}dx^{(1)}\]
Now we have got Claim 3. 

\hspace{74ex}$\square$

\vspace{1ex}

\noindent{\bf Acknowledgement}

\vspace{1ex}

The author wish to thank Joe Suzuki for reading this paper in detail and giving a lot of useful advice.

\end{document}